\documentclass{amsart}

\usepackage{amsmath}
\usepackage{amssymb}
\usepackage{mathrsfs}
\usepackage{url}
\usepackage{enumitem}
\usepackage{stmaryrd}
\usepackage[all]{xy}
\usepackage{verbatim}

\newtheorem{Thm}{Theorem}[section]
\newtheorem{Con}[Thm]{Conjecture}
\newtheorem{Cor}[Thm]{Corollary}
\newtheorem{Def}[Thm]{Definition}

\newtheorem{Eg}[Thm]{Example}
\newtheorem{Lem}[Thm]{Lemma}
\newtheorem{Prop}[Thm]{Proposition}
\newtheorem{Qn}[Thm]{Question}

\numberwithin{equation}{Thm}

\newcommand{\CC}{\mathbb{C}}

\newcommand{\FF}{\mathbb{F}}

\newcommand{\QQ}{\mathbb{Q}}

\newcommand{\ZZ}{\mathbb{Z}}
\newcommand{\PP}{\mathbb{P}}

\date{\today}

\author{Daqing Wan}
\address{Department of Mathematics, University of California, Irvine, CA 92697-3875}
\email{dwan@math.uci.edu}

\title{Zeta Functions of $\ZZ_p$-Towers of Curves}

\begin{document}

\begin{abstract}
In these notes, we explore possible stable properties for the zeta function of a geometric $\ZZ_p$-tower of curves over a finite field of characteristic $p$, 
in the spirit of Iwasawa theory. A number of fundamental questions and conjectures are proposed for those $\ZZ_p$ towers coming from algebraic geometry.  
\end{abstract}

\maketitle

%\section{Introduction}

\section{$\ZZ_p$-towers of curves}

Let $\FF_q$ be a finite field of $q$ elements with characteristic $p>0$ and let $\ZZ_p$ denote the ring of 
$p$-adic integers.  Consider a $\ZZ_p$-tower 
$$C_{\infty}: \cdots  \longrightarrow C_n \longrightarrow \cdots \longrightarrow C_1 \longrightarrow C_0$$
of smooth projective geometrically irreducible curves defined over $\FF_q$. 
The $\ZZ_p$-tower gives a continuous group isomorphism 
$$\rho: G_{\infty}:= {\text{Gal}}(C_{\infty}/C_0) \cong \ZZ_p.$$
For each integer $n\geq 0$,  reduction modulo $p^n$ gives an isomorphism 
$$ G_{n}:= {\text{Gal}}(C_{n}/C_0) \cong \ZZ/{p^n\ZZ}.$$
Let $S$ be the ramification locus of the tower, which is a subset $S$ of closet points of $C_0$. 
The tower is un-ramified on its complement $U = C_0 - S$. We shall assume that $S$ is finite.  

By class field theory, the ramification locus $S$ is non-empty, and for each non-empty finite $S$, there are uncountably many such $\ZZ_p$-towers over $\FF_q$. 
In fact, all $\ZZ_p$-towers over $C_0$ can be explicitly classified by the Artin-Schreier-Witt theory, see \cite{KW}. 
In contrast, if $\ell$ is a prime different from 
$p$, there are no $\ZZ_{\ell}$-towers over $\FF_q$.  Note that constant extensions do not give a tower in our sense 
since that would produce curves which are geometrically reducible.  Thus, the $\ZZ_p$-towers  
we consider are all geometric $\ZZ_p$-towers. 

The most important class of $\ZZ_p$-towers 
naturally comes from algebraic geometry. 
We describe this briefly now. 
Let $G_U$ be the Galois group of the maximal abelian extension of the function field of $U$ which is unramified on $U$. 
Let $\phi: G_U \longrightarrow {\ZZ_p^*}$ be a continuous rank one surjective $p$-adic representation. 
Let 
$$({\ZZ_p^*})^{p^n} = \{a^{p^n} | a  \in  \ZZ_p^*\}.$$
This is a subgroup of  ${\ZZ_p^*}$. The composition of $\phi$ and the reduction homomorphism give a surjective 
homomorphism (if $p>2$) 
$$\Phi: G_U \longrightarrow {\ZZ_p^*} \longrightarrow {\ZZ_p^*}/({\ZZ_p^*})^{p^n}  \cong \ZZ_p/p^n\ZZ_p.$$
%such that the composed homomorphism 
%$$\log_{(1+{\bf p})} \circ \phi: G_U { \buildrel \phi \over \longrightarrow} ~ {\ZZ_p^*} {\buildrel \log_{(1+{\bf p})} \over \longrightarrow}  \ZZ_p$$
%is surjective, where ${\bf p} =p$ if $p>2$ and ${\bf p} =4$ if $p=2$, and $\log_{(1+{\bf p})}$ denotes the $p$-adic logarithm.  
These surjective homomorphisms naturally produce a $\ZZ_p$-tower of $C_0$, unramified on $U$, which we assume to be geometric in our sense, that is, there is no constant 
subextensions. 
%Without taking the $p$-part, one would get a ${\ZZ_p^*}$-tower. 
%The most important part of $\phi$ is the $p$-part, which will be our principal concern.  
This geometric tower is further called arising (or coming) from {\bf algebraic geometry} if $\phi$ arises 
from a relative $p$-adic \`etale cohomology of an ordinary family of smooth proper variety $X$ parameterized by $U$, or 
more generally in the sense of Dwork's unit root conjecture \cite{Dw} as proved in \cite{wan}, that is,  $\phi$ comes from the unit root part of an ordinary 
overconvergent $F$-crystal on $U$.  A classical example is the Igusa ${\ZZ_p}$-tower 
arising from the universal family of ordinary elliptic curves over $\FF_p$, more generally, the ${\ZZ_p}$-tower arising from the following Dwork family of 
Calabi-Yau hypersurfaces  over $\FF_p$ parametrized 
by $\lambda$: 
$$X_0^{n+1} + X_1^{n+1}+\cdots + X_n^{n+1} + \lambda X_0X_1\cdots X_n =0.$$

\section{Zeta functions} 

For integer $n\geq 0$, let $Z(C_n, s)$ denote the zeta function of $C_n$. It is defined by 
$$Z(C_n,s) = \prod_{x\in |C_n|} \frac{1}{1-s^{{\rm deg}(x)}} \in 1 +s\ZZ [[s]],$$
where $|C_n|$ denotes the set of closed points of $C_n$. 
The Riemann-Roch theorem implies that the zeta function is a rational function in $s$ of the form 
$$Z(C_n, s) = \frac{P(C_n, s)}{(1-s)(1-qs)}, \ P(C_n, s) \in 1 +s\ZZ[s],$$
where $P(C_n, s)$ is a polynomial of degree $2g_n$ and $g_n= g(C_n)$ denotes the 
genus of the curve $C_n$. 
By the celebrated theorem of Weil, the polynomial $P(C_n, s)$ is 
pure of $q$-weight $1$, that is, the reciprocal roots of $Z(C_n, s)$ all have complex absolute value equal to $\sqrt{q}$. 

The $q$-adic valuations (also called $q$-slopes) of the reciprocal roots of $Z(C_n, s)$ remain quite mysterious in general. 
Our aim is to study this question in the spirit of Iwasawa theory. Namely, we want to investigate 
possible stable properties for the slopes when $n$ varies.  
Classical geometric Iwasawa theory corresponds to the study of the slope zero part of the zeta function. Although much more complicated, we feel that 
there is a rich theory in all higher slopes, for $\ZZ_p$-towers arising from algebraic geometry.  
The aim of these notes is to explore some of the basic questions, suggesting what might be true and 
what might be false, and giving supporting examples when available. 

\section{Main questions and conjectures}

We now make our questions a little more precise, beginning with some simpler questions. Write 
$$P(C_n, s) = \prod_{i=1}^{2g_n} (1-\alpha_i(n)s) \in \CC_p [s], \ \  0\leq v_q(\alpha_1(n)) \leq \cdots \leq v_q(\alpha_{2g_n}(n)) \leq 1, $$
where the $q$-adic valuation is normalized such that $v_q(q)=1$. We would like to understand how the polynomial $P(C_n, s)$ (and its zeros) varies when $n$ 
goes to infinity. The first and simplest question is about the degree of the polynomial $P(C_n, s)$.

\begin{Qn} How the genus $g_n$ varies 
as $n$ goes to infinity? 
\end{Qn}

The answer depends on how the $\ZZ_p$-tower is given. Classically, the construction of all $\ZZ_p$-towers was given by Witt using Witt vectors, 
and the genus was explicitly computed by Schmid (1937) in the same framework of Witt vectors, see \cite{KW} for a simplified complete  treatment.  
In modern applications, a $\ZZ_p$-tower naturally arises from algebraic geometry, and 
it is not clear how to compute the genus.  In this direction, we propose 

\begin{Con}\label{Con1.0}
Assume that the tower comes from  algebraic geometric. Then, the genus sequence $g_n$ is stable, that is, 
there are constants $a, b, c$ with $a>0$ depending on the tower such that for all sufficiently large $n$, we have 
$$g_n = ap^{2n} + bp^n + c.$$
\end{Con}

As an example, one can check that the conjecture is true in the case of Igusa $\ZZ_p$-tower, using the results in Katz-Mazur [KM].  
The converse of the conjecture is not true. There are only countably many $\ZZ_p$-towers coming from algebraic geometry, but there 
are uncountably many $\ZZ_p$-towers with stable genus sequence. 
The genus growth question was first studied in Gold-Kisilevsky \cite{GK}, where a simple lower bound of 
the form $g_n \geq cp^{2n}$ is proved for some positive constant $c$.   
It is easy to construct towers so that $g_n$ grows as fast as one wishes. In \cite{KW}, we give 
a general explicit formula for $g_n$ via a different explicit construction of $\ZZ_p$-towers. This new construction is better suited for studying   
the slopes of zeta functions. In this way,  one can quickly read off the genus and see its variational behavior.  In particular, genus stable 
$\ZZ_p$-towers are completely classified in \cite{KW} \cite{KW2} in terms of a simple condition on our explicit Witt construction of $\ZZ_p$-towers. To prove the above 
conjecture, one needs to check that this simple condition is satisfied for all towers coming from algebraic geometry, which are usually in multiplicative forms. 
This is not easy to do in practice. In \cite{K2}, Joe Krammer-Miller develops  a multiplicative approach which is better suited to prove this conjecture. 

{\bf Remarks}. One can ask the same question for the stable behavior of  the discriminant of a $\ZZ_p$-tower of number fields. Let 
$$K= K_0 \subset K_1 \subset K_2 \subset \cdots $$ be a $\ZZ_p$-tower of number fields. Let $d_n$ denote the absolute discriminant 
of $K_n$. The analogue of the genus in a number field is the logarithm of the discriminant. 
Since the prime to $p$-part of $d_n$ is fixed for large $n$, this is equivalent to asking for the $p$-adic valuation $v_p(d_n)$. 
One can show \cite{Tate} that for 
sufficiently large $n$, 
$$v_p (d_n)  = anp^n +bp^n +c, $$
where $a>0, b, c$ are constants depending only on the tower.  
In the beautiful paper \cite{Up}, Upton has generalized this stable formula  to any $p$-adic Lie extension of number fields, not necessarily 
abelian. For the discriminant problem, the number field case is easier as the ramification is much more controlled. 
In contrast,  the function field case is much more complicated as the ramification is completely wild, and so some condition is necessary to 
insure stability.

We now return to the function field case. 
Having understood the degree of the polynomial $P(C_n, s)$, our second natural question is about the splitting field of $P(C_n, s)$. 
In particular, how the extension degree of the splitting field of $P(C_n, s)$ varies as $n$ varies. 
\begin{Con}\label{Con1.1}
Let $\QQ_n$ (resp. $\QQ_{p,n}$) denote the splitting field of $P(C_n, s)$ over $\QQ$ (resp. $\QQ_p$). 

\item{(1)}. The extension degree $[\QQ_{n}:\QQ]$ goes to infinity  
as $n$ goes to infinity. 

\item{(2)}. The extension degree $[\QQ_{p,n}:\QQ_p]$ goes to infinity 
as $n$ goes to infinity. 

\item{(3)}. The ramification degree $[\QQ_{p,n}:\QQ_p]^{\rm{ram}}$ goes to infinity 
as $n$ goes to infinity. 

\item{(4)}. There is a positive constant $c$ depending on the tower such that for all large $n$, we have 
$[\QQ_{p,n}:\QQ_p]^{\rm{ram}}\geq cp^n$. 
\end{Con}

Clearly, each part is significantly stronger than its previous part.  In all the examples we know, 
the ramification degree of $\QQ_{p,n}$ over $\QQ_p$ is at least as large as $cp^n$ as $n$ grows, that is, the strongest part holds 
and thus all parts hold. Note that this conjecture is for all $\ZZ_p$-towers, not 
necessarily those coming from algebraic geometry. If the $\ZZ_p$-tower comes from algebraic geometry, then 
one expects the much stronger stability that $[\QQ_{p,n}:\QQ_p]^{\rm{ram}}= cp^n$ for some positive constant $c$ and all large $n$. 
This would be a consequence of our later slope stability conjecture. The work of \cite{GK} shows that the 
above conjecture is true if $\mu=0$ and $\lambda\not=0$, where $\mu$ and $\lambda$ are the invariants in Iwasawa theory 
which will be discussed later. 
%The above conjecture can be shown to be true if $\lambda\not=0$ ($\mu$ can be arbitrary). 

We now refine the above geometric genus question into many arithmetic questions in terms of slopes. 
Fix a rational number $\alpha\in [0, \infty )$, let $d_{\alpha}(n)$ 
denote the multiplicity of $\alpha/p^{n}$ in the slope sequence of $P(C_n, s)$. That is, 
$$d_{\alpha}(n) = \# \{ 1\leq i \leq 2g_n | v_q(\alpha_i(n))=\frac{\alpha}{p^{n}}\}.$$
The reason to re-scale the slope by the factor $1/p^{n}$ is that  the first few slopes of $P(C_n, s)$ is expected to be of the form $\alpha/p^n$ for a few rational 
numbers $\alpha$, independent of $n$, if the tower comes from algebraic geometry.   By definition, $d_{\alpha}(n)=0$ for $\alpha>p^n$, and 
$$\sum_{\alpha} d_{\alpha}(n) = 2g_n.$$

\begin{Qn} For a fixed rational number $\alpha \in [0,\infty )$, how the number $d_{\alpha}(n)$ varies 
 as $n$ goes to infinity? 
\end{Qn}

Note that $d_0(n)$ is the $p$-rank of of $C_n$, namely, the rank of the $p$-adic Tate module of the Jacobian of $C_n$. 
Furthermore, the class number $h_n$ (the number of divisor classes 
of degree $0$ of $C_n$ ) is given by the residue formula 
$$h_n: = \#{\rm Jac}(C_n)(\FF_q) = P(C_n, 1).$$
One has the following result from 
classical geometric Iwasawa theory. 

\begin{Thm}\label{Thm1.1}  (i). There are integer constants $\mu_1$ and $\mu_2$ such that for all sufficiently large $n$, we have 
$$d_0(n) = p^n\mu_1 + \mu_2.$$
 (ii). There are integer constants $\mu, \lambda, \nu$ such that for all sufficiently large $n$, we have 
$$v_p(h_n) = \mu p^n + \lambda n +\nu.$$
(iii). For all sufficiently large $n$, we have the congruence 
$$\frac{h_n/h_{n-1}}{p^{v_p(h_n/h_{n-1})}} \equiv 1 \mod p.$$
\end{Thm}

This theorem is really about the slope zero part of the zeta function. To illustrate our new point of views, 
we shall give a simple new proof of the above result, using the definition of the $T$-adic L-function introduced later. 
This new approach is successfully extended in \cite{W19} to prove the functional field version of Greenberg's conjecture for 
$\ZZ_p^r$-towers for any positive integer $r$. 
It is a pleasure to thank Ralph Greenberg who asked if part (i) of the above theorem is true. 
For higher slope $\alpha>0$, the problem is completely new and apparently more complicated.

\begin{Qn}\label{Qn1.3}
For each fixed $\alpha \in [0,\infty )$, are there constants $\mu_1(\alpha)$ and $\mu_2(\alpha)$ such that for all sufficiently large $n$, we have 
$$d_{\alpha}(n) = p^n\mu_1(\alpha) + \mu_2(\alpha)?$$
\end{Qn}

If we do not re-scale the slope in the definition of $d_{\alpha}(n)$ and look at all slopes, not necessarily the first few slopes, we can ask the following question. 
%Yet, we expect that a good understanding of 
%higher slope part of the zeta function would shed new light 
%on the variation of the special values of $P(C_n,s)$ at $s=1$ and at many other points, when $n$ varies.  Iwasawa's original motivation 
%comes from much simpler $\ZZ_p$-towers obtained by constant extensions, and hence not a geometric tower in our sense. 

\begin{Qn} Are the $q$-slopes  
$$\{ v_q(\alpha_1(n)), \cdots, v_q(\alpha_{2g_n}(n))\} \subset [0, 1] \cap \QQ \subset [0, 1]$$ 
equi-distributed in the interval $[0, 1]$ as $n$ goes to infinity? 
\end{Qn}

A weaker version is to ask if the set of $q$-slopes for all $n$ is dense in $[0, 1]$.  
A stronger version is the following possible finiteness property. 

\begin{Qn} Is the slope sequence stable in some sense? More precisely, is there a positive integer $n_0$ depending on the tower such that 
the re-scaled $q$-slopes
$$\{ p^nv_q(\alpha_1(n)), \cdots, p^nv_q(\alpha_{2g_n}(n))\}$$
for all $n> n_0$ are determined explicitly in a simple way by 
their values for $0\leq n \leq n_0$, using a finite number of arithmetic progressions? 
\end{Qn}

The precise meaning of this finiteness property  will be made clearer later. 
These three questions are too general to have a positive answer  in full generality, as there are too many 
$\ZZ_p$-towers, most of them are not natural.  In fact, we believe that each of the above three questions has a 
negative answer in full generality.  It would be interesting to find examples showing that 
the above three questions indeed have a negative answer.  
However, we conjecture that the answers to all three  questions are positive for all  
$\ZZ_p$ towers coming from algebraic geometry. 
More precisely, our main conjecture is the following. 

\begin{Con}\label{Con1.2}
Assume that the tower comes from algebraic geometry.

\item{(1)}. The genus sequence $g_n$ is stable. That is, 
there are constants $a, b, c$ with $a>0$ depending on the tower such that for all sufficiently large $n$, we have 
$$g_n = ap^{2n} + bp^n + c.$$

\item{(2)}. For each fixed $\alpha \in [0,\infty )$, there are constants $\mu_1(\alpha)$ and $\mu_2(\alpha)$ such that for all sufficiently large $n$, we have 
$$d_{\alpha}(n) = p^n\mu_1(\alpha) + \mu_2(\alpha).$$

\item{(3)}. The $q$-slopes  
$$\{ v_q(\alpha_1(n)), \cdots, v_q(\alpha_{2g_n}(n))\} \subset [0, 1] \cap \QQ \subset [0, 1]$$ 
are equi-distributed in the interval $[0, 1]$ as $n$ goes to infinity. 

\item{(4)}. The slope sequence is stable. That is, there is a positive integer $n_0$ depending on the tower such that the re-scaled $q$-slopes
$\{ p^nv_q(\alpha_1(n)), \cdots, p^nv_q(\alpha_{2g_n}(n))\}$ for all $n> n_0$ are determined explicitly by 
their values for $0\leq n \leq n_0$, using a finite number of arithmetic progressions. 
\end{Con}

We shall make part (4) of the conjecture more precise later and give some supporting examples.  
Parts (2)-(4) of the conjecture are in increasing level of difficulties. 
%In particular, we expect that part (1) of the conjecture would be provable within foreseeable future. 
Each of parts (1)-(3) is a weak 
consequence of part (4). Thus, a proof for each of parts (1)-(3) is also an evidence for the strongest part (4) of 
the conjecture.  Recent work of Kosters-Zhu \cite{KZ} suggests that part (1) almost implies part (3), which has been 
proved when $U$ is the affine line. 
%Note that part (4) already implies a stronger form of Conjecture \ref{Con1.1} for $\ZZ_p$ towers coming from algebraic geometry. 

Before moving on, we note the following consecutive divisibility 
$$P(C_0, s) | P(C_1, s) | \cdots | P(C_n, s) | \cdots.$$
Thus, it is enough to study for $n\geq 1$, the primitive part of $Z(C_n, s)$ defined by 
$$Q(C_n, s) = \frac{P(C_n, s)}{P(C_{n-1}, s)} = \frac{Z(C_n, s)}{Z(C_{n-1}, s)}.$$

\section{L-functions} 

Note that if $C_n$ is ramified over $C_0$ at a closed point $x\in S$, then $C_{m}$ is totally ramified over $C_n$ at $x$ 
for all $m\geq n$ since the Galois group is a cyclic $p$-group.  Without loss of generality, by going to a larger $n$ if necessary, 
we can assume that $C_1$ is already ramified at every point of $S$. 
From now on, we assume that $C_1$ is indeed (totally) ramified at every point of $S$. 

Recall that for $n\geq 1$, the Galois group 
$$G_{n}= {\text{Gal}}(C_{n}/C_0) \cong \ZZ/{p^n\ZZ}.$$
For a primitive character $\chi_n: G_n \rightarrow \CC_p^*$ of order $p^n>1$, the L-function of $\chi_n$ over $C_0$ is 
$$L(\chi_n, s) = \prod_{x\in |U|} \frac{1}{1-\chi_n(\text{Frob}_x)s^{\text{deg}(x)}} \in 1 +s\ZZ[\chi_n(1)][[s]],$$
where $|U|$ denotes the set of closed points of $U$ and $\text{Frob}_x$ denotes the arithmetic Frobenius 
element of $G_n$ at $x$.  Note that $\zeta_{p^n}:=\chi_n(1)$ is a primitive $p^n$-th root of unity. 
Again, Weil's theorem shows that the L-function $L(\chi_n, s)$ is a polynomial 
in $s$, pure of weight $1$. 
One has the decomposition 
$$Q(C_n, s) = \prod_{\chi_n: G_n \rightarrow \CC_p^*} L(\chi_n, s),$$
where $\chi_n$ denotes a primitive character of order $p^n$.  
For $\sigma \in \text{Gal}(\QQ(\zeta_{p^n})/\QQ) = \text{Gal}(\QQ_p(\zeta_{p^n})/\QQ_p)$, one checks that 
$$L(\chi_n, s)^{\sigma}  = L(\chi_n^{\sigma}, s).$$
It follows that the degree and the slopes for $L(\chi_n, s)$ depend only on $n$, not on the choice of 
the primitive character $\chi_n$ of $G_n$.  We can just choose and fix one character $\chi_n$ of order $p^n$ 
for each $n\geq 1$, if desired.  

The extension degree conjecture for the splitting field of $P(C_n, s)$ is reduced to 

\begin{Con}\label{Con2}

Let $L_{p,n}$ denote the splitting field of $L(\chi_n, s)$ over $\QQ_p$.

\item{(1)}. The extension degree $[L_{p,n}:\QQ_p]$ goes to infinity 
as $n$ goes to infinity. 

\item{(2)}. The ramification degree $[L_{p,n}:\QQ_p]^{\rm{ram}}$ goes to infinity 
as $n$ goes to infinity. 

\item{(3)}. There is a positive constant $c$ depending on the tower such that for all large $n$, we have 
$[L_{p,n}:\QQ_p]^{\rm{ram}}\geq cp^n$. 
\end{Con}

%As noted before, this conjecture is now known to be true if the $\lambda$-invariant of the tower is 
%non-zero. 

Let $\ell(n)$ denote the degree of the polynomial $L(\chi_n, s)$. The degree of $Q(C_n, s)$ is simply 
$p^{n-1}(p-1)\ell(n)$. The genus $g_n$ of $C_n$ is given by the formula 
$$2g_n -2g_0= (p-1)\sum_{i=1}^n p^{i-1}\ell(i). $$
To understand the variation of $g_n$ as $n$ varies, it is enough to understand the 
variation of $\ell(n)$ as $n$ varies. Conjecture \ref{Con1.0} is equivalent to 

\begin{Con}\label{Con2.0}
Assume that the tower comes from  algebraic geometry. Then, the degree sequence $\ell(n)$ is stable in the sense that 
there are constants $a, b$ with $a>0$ depending on the tower such that for all sufficiently large $n$, we have 
$$\ell(n) = ap^{n} + b.$$
\end{Con}

Again, this conjecture can be refined in terms of slopes. 
Write 
$$L(\chi_n, s) = \prod_{i=1}^{\ell(n)} (1-\beta_i(n)s) \in \CC_p [s], \ \  0\leq v_q(\beta_1(n)) \leq \cdots \leq v_q(\beta_{\ell(n)}(n)) \leq 1.$$
The slope sequence for $Q(C_n, s)$ is just equal to the slope sequence for $L(\chi_n, s)$, repeated $p^{n-1}(p-1)$ times. 

Fix a rational number $\alpha\in [0, \infty )$, let $\ell_{\alpha}(n)$ 
denote the multiplicity of $\alpha/p^n$ in the slope sequence of $L(\chi_n, s)$. That is, 
$$\ell_{\alpha}(n) = \# \{ 1\leq i \leq \ell(n) | v_q(\beta_i(n))=\frac{\alpha}{p^n}\}.$$
It is clear that for every $n\geq 1$, we have 
$$\sum_{\alpha} \ell_{\alpha}(n) = \ell(n),$$
and for every $\alpha$, we have 
$$(p-1)\sum_{i=1}^n p^{i-1}\ell_{\alpha}(i) = d_{\alpha}(n) - d_{\alpha}(0). $$

\begin{Qn} For a fixed rational number $\alpha \in [0,\infty )$, how the number $\ell_{\alpha}(n)$ varies 
 as $n$ goes to infinity? 
\end{Qn}

In the case $\alpha = 0$, the situation is quite simple and clean. 
The following result follows quickly from the definition of $T$-adic L-functions introduced later. 

\begin{Lem}\label{Lem1.1} For every $n\geq 1$, we have $\ell_0(n) =d_0(0)-1 + \sum_{x\in S} {{\rm deg}(x)}$. 
\end{Lem}

In particular, the number $\ell_0(n)$ for $n\geq 1$ is a constant independent of $n$. It follows that 
$$d_0(n) - d_0(0) = \ell_0(1)\sum_{i=1}^n (p-1)p^{i-1} =  \frac{p^n-1}{p-1}(d_0(1)-d_0(0)).$$
Substituting the value of $\ell_0(1)$ in the lemma, we get an alternative formula for $d_0(n)$: 
$$d_0(n) -d_0(0) = (p^n-1)(d_0(0) -1 + \sum_{x\in S} {{\rm deg}(x)}).$$
This proves part (i) of Theorem \ref{Thm1.1}. 
It implies that the sequence $\{d_0(n)\}_{n\geq 0}$ is determined by the first two terms $d_0(0)$ and $d_0(1)$. 
Alternatively, it is determined by the first term $d_0(0)$ and the degree of the divisor $S$. 
For $\alpha>0$, the slope sequence $\{d_{\alpha}(n)\}_{n\geq 0}$  is more complicated. 
Similar to the zeta function case, we can ask 

\begin{Qn}\label{Qn2.1}

\item{(1)}. 
For each fixed $\alpha \in [0,\infty )$, the number $\ell_{\alpha}(n)$ is a constant for all sufficiently large $n$?  

\item{(2)}. As $n$ goes to infinity, are the $q$-slopes 
$$\{ v_q(\beta_1(n)), \cdots, v_q(\beta_{\ell(n)}(n))\} \subset [0, 1] \cap \QQ \subset [0, 1]$$ 
equi-distributed in the interval $[0, 1]$?  

\item{(3)}. 
Is there a positive integer $n_0$ depending on the tower such that the re-scaled $q$-slopes  
$\{ p^nv_q(\beta_1(n)), \cdots, p^nv_q(\beta_{\ell(n)}(n))\}$ for all $n> n_0$ are determined explicitly by 
their values for $0\leq n \leq n_0$, using a finite number of arithmetic progressions? 
\end{Qn}

Again, these questions are too general to have a positive answer in full generality. 
However, we conjecture all of them have a positive answer for $\ZZ_p$-towers coming from algebraic geometry.  
To gain some feeling about what part (3) means, we give an example next.

\section{$\ZZ_p$-towers over the affine line} 
In this section, we explain that all the above questions for all slopes have a simple positive answer 
for many $\ZZ_p$-towers over the affine line, as first studied in \cite{DWX}, and more recently in \cite{KZ}\cite{Li}.  

Fix an element $\beta$ of $W(\FF_q)=\ZZ_q$ with trace $1$. By the classification in \cite{KW},  any $\ZZ_p$-tower over $C_0=\PP^1$, 
totally ramified at infinity $\infty$ and unramified on $\mathbb{A}^1 =\PP^1 -\{ \infty\}$ can be uniquely constructed from a constant $c\in \ZZ_p$ and 
a primitive convergent power series 
$$f(x) = \sum_{(i,p)=1} c_i x^i \in \ZZ_q[[x]], \ c_i \in \ZZ_q, \  \lim_i c_i =0,$$ 
where $f(x)$ is called primitive if not all $c_i$ are divisible by $p$, that is, $f(x)$ is not divisible by $p$. 
The construction is explicitly given by the following equation 
$$C_{\infty}: [y_0^p, y_1^p, \cdots] - [y_0, y_1, \cdots] = c\beta + \sum_{(i,p)=1}c_i [x^i, 0, \cdots],$$
where both sides are Witt vectors. The constant $c$ has no contribution to all our questions, and thus we shall asumme that $c=0$. 
The tower is then uniquely constructed from the primitive convergent power series $f(x)$. 

Comparing the coordinates, one finds that $C_n$ is defined by a system of $n$ polynomial equations over $\FF_q$ in $n+1$ variables $(y_1,\cdots, y_n, x)$. 
It is clear that $C_0 = \PP^1$, $C_1$ is the 
Artin-Schreier curve whose affine equation over $\FF_q$ is given by 
$$y_0^p -y_0 = f(x),$$
and $C_2$ is the curve above $C_1$ given by an additional equation (over $\FF_q$)
\[
y_1^p - y_1 + \frac{y_0^{p^2}-y_0^p - (y_0^p - y_0)^p}{p} =  \frac{f^\sigma(x^p) - f(x)^p}{p},
\]
where $ f^\sigma(x) = \sum_{(i,p)=1} \sigma(c_i) x^i$ and $\sigma$ is the Frobenius automorphism on $\ZZ_q$. 
Note that over $\FF_q$, the power series $f(x)$ reduces to a polynomial. 

The map $C_n \longrightarrow C_{n-1}$ is given by the projection 
$$(y_1, \cdots, y_n, x) \longrightarrow (y_1, \cdots, y_{n-1}, x).$$
Let $\zeta_{p^n}$ be a primitive $p^n$-th root of unity in $\CC_p^*$. Set $t_n = \zeta_{p^n}-1$, which is a uniformizer in the local field 
$\QQ_p(\zeta_{p^n})$. 
For all $n\geq 1$, it is clear that $\chi_n \equiv 1 \mod (t_n)$ and 
$$L(\chi_n, 1) \equiv Z(\mathbb{A}^1, 1) \equiv 1 \mod t_n.$$
We obtain 
\begin{Thm}  Let $h_n$ denote the class number of $C_n$. Then, for all $n\geq 0$, we have 
$$h_n \equiv 1 \mod (p).$$
In particular, $d_0(n) = \ell_0(n) = \mu = \lambda = \nu =0$. 
\end{Thm} 
It would be interesting to explore possible improvements of the above congruence to a congruence 
modulo higher powers of $p$, that is, to understand the first few digits in the $p$-adic expansion of $h_n$. 

For integer $n\geq 1$, the Artin conductor $a(\chi_n)$ of the character $\chi_n$ is calculated explicitly in \cite{KW}: 
$$a(\chi_n) = 1 + \max_{v_p(c_i)<n} \{ ip^{n-1-v_p(c_i)}\}.$$
It follows that the degree  of $L(\chi_n, s)$ is given by 
the formula
$$\ell(n) = -1 + p^{n-1} \max_{v_p(c_i)<n} \{ ip^{-v_p(c_i)}\}.$$
This tower is genus stable if and only if $\ell(n)$ is a linear polynomial in $p^n$ for large $n$. 
This is the case when 
$$d: = \max_{(i,p)=1} \{ \frac{i}{p^{v_p(c_i)}}\}$$ exits and is a finite rational number, in which case, for all sufficiently large $n$, we have the stable degree formula 
$$\ell(n) = dp^{n-1} -1.$$
This case is called {\bf strongly genus stable}. There is a more complicated class of genus stable towers \cite{KW2} that we do not discuss here for 
simplicity.

\begin{Def} Assume that the tower is strongly genus stable as above. Let $n_0$ be the smallest positive integer $k$ such that 
$\ell(k)=p^{k-1}d-1$, and if 
 the $q$-slope sequence of $L(\chi_k, s)$ is denoted by $\{\alpha_1, \cdots, \alpha_{p^{k-1}d-1}\}$,  
then for every integer $n\geq k$, the $q$-slope sequence of $L(\chi_n, s)$ is given by 
\[
\bigcup_{i=0}^{p^{n-k}-1} \big\{\frac{i}{p^{n-k}}, \frac{\alpha_1+i}{p^{n-k}}, 
\dots, \frac{\alpha_{dp^{k-1}-1}+i}{p^{n-k}}\big\} - \{0\}. 
\]
If such positive integer $k$ does not exist, we define $n_0 =\infty$. 
\end{Def}

The finiteness of the number $n_0$ hence implies a striking stable property for the slope sequence of 
$L(\chi_n, s)$ as $n$ goes to infinity. It implies that the slopes of $L(\chi_n, s)$ normalized by the factor $p^n$ for all large $n$ are given by 
a finite number of arithmetic progressions. For this reason, the tower is called {\bf slope stable} if $n_0$ is finite. 
Note that if the tower is slope stable, then clearly the tower must be genus stable. It is tempting to  conjecture that the converse is also true. 
Although we do not have a counter-example, 
we are a little cautions here and  
will just state it as a question.

\begin{Qn} Assume the tower is genus stable. 
	Is the tower slope stable? 
\end{Qn}

Note that we conjectured that the 
answer is positive for towers coming from 
algebraic geometry.  An important example is the Igusa $\ZZ_p$-tower which is genus stable, but the 
finiteness property of $n_0$ for Igusa tower seems open.  This might be related to the geometry of eigencurves in the framework of modular forms, 
see \cite{wan-xiao-zhang} \cite{liu-wan-xiao} for a recent progress 
and \cite{BP} for a precise conjectural  slope description in the case of regular primes $p$. Another related slope problem in the setting of symmetric powers of Kloosterman sums is recently studied in Haessig \cite{Ha}, and Fresan-Sabbah-Yu \cite{FSY}. 
We now give some examples of $\ZZ_p$-towers constructed using the primitive convergent power series $f$, where 
the above question has a positive answer. 

\begin{Def} 
\item{(1)}. The tower is called a {\bf polynomial tower} of degree $d$ if the primitive convergent power series $f(x)=\sum_{(i,p)=1} c_i x^i$ is a polynomial of degree $d$. 

\item{(2)}. 
The tower is called a {\bf unit root tower} of degree $d$ if $f(x)$ is a polynomial of degree $d$ and furthermore all non-zero 
coefficients $c_i$ are roots of unity.  
\end{Def} 

Clearly, a polynomial tower is strongly genus stable and its degree $d$ is not divisible by $p$. 

\begin{Thm}[\cite{DWX}]  For a unit root tower of degree $d$, the number $n_0$ is finite. Furthermore, 
we have the following explicit upper bound 
$$n_0 \leq 1 + \lceil \log_p(\frac{d}{8} \log_pq)\rceil.$$
In particular, if $p\geq \frac{d}{8} \log_pq$, then $n_0 \leq 2$. 
\end{Thm} 

\begin{Cor} For a unit root tower of any degree, we have the following 

\item{(1)}. The $q$-slopes of $L(\chi_n, s)$ (resp., $P(C_n, s)$) are equi-distributed in $[0, 1]$ as $n$ goes to infinity. 

\item{(2)}. For every rational number $\alpha \in [0, \infty )$, the sequence $\ell_{\alpha}(n)$ is a constant independent of 
$n \geq n_0$. 

\item{(3)}. The ramification degree of $P(C_n, s)$ over $\QQ_p$ is equal to $cp^n$ for some positive constant $c$ for all $n>n_0$.  

\end{Cor}

The above explicit bound for $n_0$ can be improved in various special cases. 

\begin{Eg} For a unit root tower of degree $d$ satisfying $p\equiv 1 \mod d$, we have $n_0=1$, $\alpha_i = i/d$ for $1\leq i \leq d$ and hence the slope sequence 
of $L(\chi_n, s)$ for all $n\geq 1$ is given by 
$$\{ \frac{1}{dp^{n-1}},   \frac{2}{dp^{n-1}}, \cdots,  \frac{dp^{n-1}-1}{dp^{n-1}}\}. $$
This was first proved in \cite{LW}. The $T$-adic L-function introduced there plays a crucial role 
in the proof of the above more general theorem. 
\end{Eg}

\begin{Eg} Let $f(x) = x^d +ax \in \ZZ_q[x]$ with $a^{q-1} =1$, $d$ not divisible by $p$ and 
$$p > \max \{ \frac{d^2(d-1)}{2}, 1 + \frac{d(d-1)}{4}\log_pq\}.$$
It is proved in \cite{OY}  that $n_0=1$ and 
$$\alpha_i = \frac{i}{d} + \frac{d-1}{d(p-1)} (ip - [\frac{ip}{d}]d -i).$$
\end{Eg}

\begin{Eg} If $p> 3d$, one can show that $n_0=1$ for a generic $\bar{f}$ over $\bar{\FF}_p$, see \cite{LLN}. 
\end{Eg}

The slope stable property is proved to be true for any polynomial tower in \cite{Li} and more generally for a much larger class of 
strongly genus stable towers over the affine line in [KZ].  The $T$-adic L-function introduced in next section played an essential role in  [DWX] and [Li]. 
A novel feature of [KZ] is the introduction of the $\pi$-adic L-function in infinitely many 
variables which refines and interpolates the $T$-adic L-function. 
It would be very interesting to prove the slope stable property for all strongly genus stable towers, or all genus stable towers or to find a counter-example.

\section{$T$-adic L-functions} 

We now return to the situation of an arbitrary  $\ZZ_p$-tower 
and define the $T$-adic L-function first introduced in [LW]. 
Instead of just finite characters $\chi_n: G_n \longrightarrow \CC_p^*$, we will also consider all continuous $p$-adic characters $\chi: G_{\infty} \longrightarrow \CC_p^*$, not necessarily of finite order. 
The isomorphism 
$$\rho: G_{\infty} \cong \ZZ_p$$ 
is crucial for us. The $p$-adic valued Frobenius function it induces 
$$F_{\rho}:  |U| \longrightarrow \ZZ_p, \ F_{\rho}(x) = \rho({\rm Frob}_x)$$  
determines the $\ZZ_p$-tower by class field theory. Any condition we would impose on the tower is a condition on this Frobenius function $F_{\rho}$.

Consider 
the universal continuous $T$-adic character $\ZZ_p \longrightarrow \ZZ_p[[T]]^*$ determined by sending $1$ to $1+T$. 
Composing this universal $T$-adic character of $\ZZ_p$ with the isomorphism $\rho$, we get the universal $T$-adic character of $G_{\infty}$: 
\[
\rho_T: G_{\infty} \longrightarrow \ZZ_p \longrightarrow \text{GL}_1(\ZZ_p[[T]]) = \ZZ_p[[T]]^*.
\]
Let $D_p(1)$ denote the open unit disk in $\CC_p$. 
For any element $t \in D_p(1)$, we have a natural evaluation map $\ZZ_p[[T]]^* \rightarrow \CC_p^*$ sending $T$ to $t$.  Composing all these maps, we get, for fixed $t\in D_p(1)$, a 
continuous character
\begin{equation} \label{varying character}
\rho_t: G_{\infty} \longrightarrow \CC_p^*. 
\end{equation}
The open unit disk $D_p(1)$ parametrizes all continuous $\CC_p$-valued  characters $\chi$ of $G_{\infty}$ via the relation 
$t = \chi(1)-1$. The L-function of $\rho_t$ is defined in the usual way: 
$$L({\rho_t}, s) = \prod_{x\in |U|} \frac{1}{1-\rho_t({\rm Frob}_x)s^{{\rm deg}(x)}}\in 1 +s\CC_p[[s]].$$

In the case that $\chi=\chi_n$ is a finite $p$-adic character of $G_{\infty}$ of order $p^n$, then $\chi(1)$ is a primitive $p^n$-th roots of unity and 
we have 
$$t=t_n:= \chi_n(1)-1, \ |t_n|_p  = p^{-\frac{1}{p^{n-1}(p-1)}}, \ L(\chi_n, s) = L({\rho}_{t_n}, s). $$ 
Elements of the form $t_n =\chi_n(1)-1$ for $n\geq 0$ are called the classical points in $D(1)$. 
As $n$ goes to infinity, $t_n$ approaches to the boundary of the disk $D_p(1)$. 
Thus, to understand the behavior of $L(\chi_n, s)$ as $n$ grows, its is enough to understand the 
L-function $L(\rho_t, s)$ for all $t$ near the boundary of $D_p(1)$.  More precisely, we should understand the 
following universal L-function. 

\begin{Def} The $T$-adic L-function of the tower is the L-function of the $T$-adic character $\rho_T$: 
$$L_{\rho}(T,s):= L(\rho_T, s) = \prod_{x\in |U|} \frac{1}{1-(1+T)^{\rho({\rm{Frob}}_x)}s^{{\rm{deg}}(x)}}\in 1 +s\ZZ_p[[T]][[s]].$$
\end{Def}
This is a $p$-adic power series in the two variables $T$ and $s$. 
For $t\in D_p(1)$, we have 
$$L(\rho_{t}, s) = L_{\rho}(T,s)|_{T=t} =L_{\rho}(t, s).$$
As noted above, the specialization of $L_{\rho}(T, s)$ at every classical point $T=t_n$ is a rational function 
$L_{\rho}(t_n, s)$ in $s$, in fact, a polynomial in $s$ of degree $\ell(n)$ for $n\geq 1$. In this case, the character $\rho_{t_n}$ is of finite 
order. If $t\in D_p(1)$ is not a classical point, i.e., $\rho_t$ is of infinite order, we do NOT know a single example for which $L_{\rho}(t, s)$ is rational. 

\begin{Qn}Is  there a non-classical $t\in D_p(1)$ such that $L_{\rho}(t, s)$ is rational? 
\end{Qn}

To see the significance of the $T$-adic L-function, we now use the definition  of $L_{\rho}(T, s)$ to prove Theorem \ref{Thm1.1}. 
Since the character $\rho_T$ is trivial modulo $T$, the L-function $L_{\rho}(T,s)$ modulo $T$ is the same as the zeta function $Z(U,s)$ of $U$. 
This gives the congruence 
$$L_{\rho}(T, s) \equiv Z(C_0, s) \prod_{x\in S} (1- s^{{\rm deg}(x)})    = \frac{P(C_0, s)}{(1-s)(1-qs)} \prod_{x\in S} (1- s^{{\rm deg}(x)}) \mod T .$$
Replacing $T$ by $t_n$ for $n\geq 1$, we deduce 
$$L_{\rho}(t_n, s) \equiv \frac{P(C_0, s)}{(1-s)} \prod_{x\in S} (1- s^{{\rm deg}(x)}) \mod t_n.$$ 
Comparing the number of reciprocal roots of slope zero, one finds that for $n\geq 1$, 
$$\ell_0(n) = d_0(0) -1 +\sum_{x\in S} {{\rm deg}(x)}.$$ 
This proves Lemma \ref{Lem1.1} and hence part (i) of Theorem \ref{Thm1.1}. 

To get the stable  formula for $v_p(h_n)$, we need to specialize $s$ at $1$. Write 
$$L_{\rho}(T,s) =\sum_{k=0}^{\infty} L_k(T) s^k, \ L_k(T) \in \ZZ_p[[T]].$$
Since $L_{\rho}(t_n, s)$ is a polynomial of degree $\ell(n)$, we have $L_k(t_n)=0$ 
for all $k>\ell(n)$. The $p$-adic Weierstrass preparation theorem implies that 
$$L_k(T)  = \frac{(1+T)^{p^n}-1}{T} u_k(T), \ u_k(T)\in  \ZZ_p[[T]]$$
for all $k> \ell(n)$. It follows that the series $L_{\rho}(T,s)$ is $(p, T)$-adically convergent for $s\in \ZZ_p[[T]]$. 
Taking $s=1$, noting that $L_{\rho}(T, 1) \not =0$ as its specialization at classical points $t_n$ is non-zero, we can write 
$$L_{\rho}(T, 1) = p^{\mu}(T^{\lambda}+pa_1T^{{\lambda} -1} + \cdots +pa_{\lambda}) u(T),$$
where $a_i\in \ZZ_p$ and $u(T)$ is a unit 
in $\ZZ_p[[T]]$. It follows that   
$$v_p(h_m) -v_p(h_0)= \sum_{n=1}^m (p-1)p^{n-1}v_p(L_{\rho}(t_n, 1)) .$$
Since $v_p(u(t_n))=0$, and for $p^{n-1}(p-1) > \lambda$, we have
$$v_p(t_n^{\lambda}+pa_1t_n^{{\lambda} -1} + \cdots +pa_{\lambda}) = \frac{ \lambda}{p^{n-1}(p-1)},$$
we conclude the stable formula that for $m$ sufficiently large,
$$v_p(h_m) =\mu p^m + \lambda m + \nu$$
for some constant $\nu$. Part (ii) of Theorem \ref{Thm1.1} is proved. 

Finally, note that for $n\geq 1$, we have 
$$\frac{h_n}{h_{n-1}} = {\rm Norm}_{\QQ_p(t_n)/\QQ_p}(L_{\rho}(t_n, 1)).$$
The minimal polynomial of $t_n$ is the $p$-Eisenstein polynomial 
$$ \frac{(1+T)^{p^n}-1}{(1+T)^{p^{n-1}}-1}  = T^{p^{n-1}(p-1)} + \cdots + p.$$
Thus, for $n\geq 2$ (or $p>2$), we have 
$${\rm Norm}_{\QQ_p(t_n)/\QQ_p}(t_n) = p, \ {\rm Norm}_{\QQ_p(t_n)/\QQ_p}(u(t_n)) \equiv 1 \mod p, $$
where $u(T) $ is any unit in $\ZZ_p[[T]]$. It follows that for sufficiently large $n$, 
$$\frac{h_n/h_{n-1}}{p^{v_p(h_n/h_{n-1})}} \equiv  {1\over p^{\lambda}} {\rm Norm}_{\QQ_p(t_n)/\QQ_p}((t_n^{\lambda} +pa_1t_n^{{\lambda} -1} + \cdots +pa_{\lambda})u(t_n))\equiv 1 \mod p.$$
Part (iii) of Theorem \ref{Thm1.1} is proved.

\section{$T$-adic Meromorphic continuation}

The power series ring $\ZZ_p[[T]]$ has three obvious topologies: the $p$-adic topology, the $T$-adic topology 
and the $(p, T)$-adic topology. All these topologies are useful to us. 
In this section, we will focus on the $T$-adic topology, which will be our starting point.  
Viewing $L_{\rho}(T, s)$ as a power series in $s$ with coefficients in the 
complete discrete valuation field $\QQ_p((T))$ with uniformizer $T$, we are interested in the $T$-adic meromorphic continuation. 
Clearly, $L_{\rho}(T, s)$ is $T$-adic analytic in the $T$-adic open unit disk $|s|_T <1$. 
One can prove 

\begin{Prop} There is a decomposition 
$$L_{\rho}(T, s) =\frac{D_1(T,s)}{D_2(T,s)},$$
where $D_i(T,s)\in 1+ s\ZZ_p[[T]][[s]]$ are $T$-adic analytic on $|s|_T \leq 1$ for $1\leq i\leq 2$. Furthermore, 
$D_2(T,s) \in 1 + qs\ZZ_p[[T]][[s]]$. 
\end{Prop}

\begin{Cor}  
\item{(i)}. $L_{\rho}(T,s)$ is $T$-adic meromorphic on the closed $T$-adic unit disk $|s|_T \leq 1$, i.e., well defined for $s\in \QQ_p[[T]]$.  

\item{(ii).} $L_{\rho}(T,s)$ is $(p, T)$-adic analytic on the closed $(p, T)$-adic unit disk $|s|_{(p, T)} \leq 1$, i.e., convergent 
for $s \in \ZZ_p[[T]]$. 
\end{Cor}

Part (ii) of the corollary was already proved in the previous section. Crew \cite{crew} further showed that the $(p, T)$-adic slope zero part of $L_{\rho}(T, s)$ 
has a cohomological interpretation in terms of $p$-adic \`etale cohomology. This is the main conjecture in function fields. 
Note that part (i) of the corollary is stronger. It cannot be deduced from the results in \cite{crew}.   
%on a curve. A short elementary proof is presented in the previous subsection. Part (i) of the Corollary also extends to the general higher 
%dimensional and higher rank case, which gives a little improvement of the result in Emerton-Kisin \cite{EK} proving a conjecture of Katz. 

To understand higher slopes of the $\ZZ_p$-tower, we need to study the analytic properties of the two variable function $L_{\rho}(T, s)$ 
beyond the closed unit disk $|s|_T \leq 1$. This leads to the following two questions.

\begin{Qn}For which tower, the L-function $L_{\rho}(t, s)$ is $p$-adic meromorphic in $|s|_p < \infty $ for all $t\in D_p(1)$. 
\end{Qn}

\begin{Qn}For which tower, the $T$-adic L-function $L_{\rho}(T, s)$ is {\bf integrally} $T$-adic meromorphic in $|s|_T < \infty $ in the sense that 
$$L_{\rho}(T, s) =\frac{D_1(T,s)}{D_2(T,s)},$$
where $D_i(T,s)\in 1+ s\ZZ_p[[T]][[s]]$ are $T$-adic analytic in $|s|_T < \infty$ for $1\leq i\leq 2$. 	
	
\end{Qn}

The second question is stronger than the first one, as the integrally $T$-adic meromorphic continuation of 
$L_{\rho}(T,s)$ implies the $p$-adic meromorphic  continuation of $L_{\rho}(t, s)$ for all $t\in D_p(1)$.  
There are examples where $L_{\rho}(t,s)$ is not $p$-adic meromorphic, see \cite{wan96}. Thus, some conditions are necessary to have a
positive answer.

Composing with the map $\ZZ_p \longrightarrow \ZZ_p^*$ sending $1$ to $1+{\bf p}$, the $p$-adic valued Frobenius function 
$F_{\rho}$ defines a rank one $p$-adic representation of the Galois group $G_U$, equivalently a rank one unit root $\sigma$-module $M_{{\rho}}$ on $U$, see \cite{wan}.  It is clear that $F_{\rho}$ and $M_{\rho}$ determine each other.  
Using the Monsky trace formula \cite{Mo}, one can prove 

\begin{Thm}\label{Thm2.1} Assume that $M_{\rho}$ is $\infty\log$-convergent on $U$. Then $L_{\rho}(T, s)$ is integrally $T$-adic 
meromorphic in $|s|_T < \infty$. It follows that $L_{\rho}(t, s)$ is $p$-adic meromorphic in 
$|s|_p < \infty $ for all $t\in D_p(1).$
\end{Thm}

This result is not general enough for applications, although both the unit root  tower and the polynomial tower do satisfy the condition of the above theorem. 
More generally, the method in \cite{wan} can be used to prove the following Coleman's generalization of the rank once case of Dwork's unit root conjecture.  
In fact, this was already worked out in Grosse-Kl\"onne \cite{EG}  in  the case that $M_{\rho}$ comes from a pure slope piece of 
some finite rank overconvergent $\sigma$-module on $U$. 

\begin{Thm}\label{Thm2.2} 
Assume that $M_{\rho}$ comes from a pure slope piece of some finite rank $\infty\log$-convergent $\sigma$-module on $U$. 
Then $L_{\rho}(T, s)$ is integrally $T$-adic meromorphic in $|s|_T < \infty$. It follows that $L_{\rho}(t, s)$ is 
$p$-adic meromorphic in 
$|s|_p < \infty $ for all $t\in D_p(1)$. 
\end{Thm}

In all natural applications arising from higher dimensional arithmetic geometry, the rank one $\sigma$-module $M_{\rho}$ 
satisfies the assumption of Theorem \ref{Thm2.2}, but usually not the 
assumption of Theorem \ref{Thm2.1}. 
These results make it possible to talk about the zeros and poles of the L-function $L_{\rho}(T, s)$ if $\rho$ comes from algebraic 
geometry.

We can now state the L-function version of our main  conjecture, which answers Question \ref{Qn2.1} and is equivalent to 
Conjecture \ref{Con1.2} (the zeta function version of our main conjecture).    

\begin{Con}
Assume that the $\ZZ_p$-tower comes from algebraic geometry.  

\item{(1)}. The L-degree is stable. That is, there are constants $a, b$ with $a>0$ such that for sufficiently lawge $n$, we have $\ell(n)=ap^n +b$. 

\item{(2)}. For each fixed rational number $\alpha \in [0,\infty )$, the number $\ell_{\alpha}(n)$ is a constant for all sufficiently large $n$. 

\item{(3)}. The $q$-slopes of $L(\chi_n, s)$ are equi-distributed in $[0, 1]$ as $n$ goes to infinity. 

\item{(4)}.  The slopes are stable.  That is, there is a positive integer $n_0$ depending on the tower such that the re-scaled $q$-slopes
$\{ p^n v_q(\alpha_1(n)), \cdots, p^n v_q(\alpha_{2g_n}(n))\}$ for all $n> n_0$ are determined explicitly by 
their values for $0\leq n \leq n_0$, using a finite number of arithmetic progressions. 

\end{Con}

Part (1) is just Conjecture \ref{Con2.0} stated previously, which is equivalent to Conjecture \ref{Con1.0} on the genus stable property. 
The remarkable work of Kosters-Zhu suggests that part (1) almost implies part (3) for any $\ZZ_p$-tower. In fact, they have proved this implication 
for all $\ZZ_p$-towers over the affine line. 
As seen above, this main conjecture is completely proven for many $\ZZ_p$-towers over the affine line. 
An important example to consider is the Igusa $\ZZ_p$-tower over $\FF_p$ for which the full conjecture seems still open. 
This example is important because of its possible connection to arithmetic of modular forms and Galois representations. In fact, part of our conjecture was  
motivated by Coleman-Mazur's question \cite{coleman-mazur} on the geometry of the eigencurve near the boundary of the weight disk, see the introduction 
in \cite{liu-wan-xiao} for additional information. Loosely speaking, the eigencurve (or its spectral curve) is the ``zero locus" of the two 
variable L-function $L_{\rho}(T,s)$. 
The eigencurve was introduced to extend Hida's theory from slope zero to all higher slopes. In one aspect, our general conjecture 
can be viewed as an attempt to extend geometric Iwasawa theory from slope zero to all higher slopes.

In these notes, for simplicity we only considered $\ZZ_p$ towers of curves, which is already sufficiently interesting. One can also consider various generalizations 
in a number of directions, for example, replacing $\ZZ_p$ by a more general compact $p$-adic Lie group \cite{W19}, and replacing curves by higher dimensional varieties. 
For example, the unit root tower with higher rank Galois group $\ZZ_p^r$ is considered in \cite{RWXY}, where some new difficulty already arises. The 
genus growth behavior of non-abelian $p$-adic towers is considered in Kramer-Miller \cite{K1}.

\end{document}